\input amstex

\documentstyle{amsppt}
\document
\input xypic
\pageheight{45pc}
\NoBlackBoxes

\define\Q{\Bbb Q}
\define\G{\Bbb G}
\define\Z{\Bbb Z}
\define\N{\Bbb N}
\define\F{\Bbb F}

\define\cal{\Cal}

\magnification 1200
\topmatter
\rightheadtext{Detecting linear dependence}
\title Detecting linear dependence by reduction maps
\endtitle
\author G. Banaszak, W. Gajda, P. Kraso\'n
\endauthor
\address
Department of Mathematics, Adam Mickiewicz University,
Pozna\'{n}, Poland
\endaddress
\email banaszak\@math.amu.edu.pl
\endemail
\address
current: Max Planck Intitute f{\" u}r Mathematik in Bonn, Bonn, Germany 
\endaddress
\email banaszak\@mpim-bonn.mpg.de
\endemail
\address
Department of Mathematics, Adam Mickiewicz University,
Pozna\'{n}, Poland
\endaddress
\email gajda\@math.amu.edu.pl
\endemail
\address
Department of Mathematics, Szczecin University,
Szczecin, Poland
\endaddress
\email
krason\@sus.univ.szczecin.pl
\endemail

\abstract We consider the local to global principle for detecting
linear dependence of points in groups of the Mordell-Weil type.~As
applications of our general setting we obtain corresponding
statements for Mordell-Weil groups of non{-}CM elliptic curves and
some higher dimensional abelian varieties defined over number fields,
and also for odd dimensional K-groups of number fields.
\endabstract
\endtopmatter

\subhead  1. Introduction
\endsubhead
\medskip

Let $E$ be an elliptic curve defined over a number field $F.$ By
the classical theorems of Mordell and Weil the group $E(F)$ of
$F$-rational points is finitely generated. The torsion part of
$E(F)$ can be computed quite efficiently in many cases: for
$F{=}\Q$ due to results of Lutz, Nagell and Mazur (see, for
example [G-SOT]) and for some larger $F$ due to results of
Kamienny, Merel and others. One of the basic facts which are
useful in calculations of the torsion is the reduction theorem
[Si, Prop. 3.1, p.176] which says that the reduction map: 
$$ r_v:\, E(F)\longrightarrow E_v(\kappa_v)
$$ (where $\kappa_v{=}{\cal O}_F/v$) is an injection, when
restricted to the torsion subgroup generated by elements of order 
relatively prime to $v$, for all primes $v$ of good reduction. 
On the other hand, it is
well-known that the rank of $E(F)$ is much more difficult to
compute (cf. [RS] for an up-to-date survey on the ranks of $E(F),$
when $F{=}\Q$).
\medskip

One of our goals in this paper is to show that the reduction maps
can be used to investigate the nontorsion part of the Mordell-Weil
group. ~Given a finite set of nontorsion points of $E(F)$, one can
ask whether it is possible to detect linear dependence among 
elements of this set by reductions.~Similar questions are of
interest for the Mordell-Weil group of any higher dimensional
abelian variety over $F,$ and also more generally - if we are
given a finitely generated abelian group $B(F)$ together with maps
$r_v: B(F)\longrightarrow B_v(\kappa_v)$ to finite groups - one
map for every finite prime of $F.$ A good example of this
situation is provided by Quillen's K-groups $K_{2n{+}1}(F),$ where
$n\geq 0,$ and the homomorphisms $$ r_v:\,
K_{2n{+}1}(F)\longrightarrow K_{2n{+}1}(\kappa_v)$$ induced by
projections ${\cal O}_F\longrightarrow \kappa_v.$ Note that
$K_{2n{+}1}(F){=}K_{2n{+}1}({\cal O}_F),$ for $n>0$ and the
K-group is finitely generated cf. [Q1]. Its rank equals the order
of vanishing of the Dedekind zeta function of $F$ at $s=-n$ cf.
[B]. In this case the target of the reduction map $r_v$ is the
cyclic group $\Z/(q_v^{n{+}1}{-}1)$, where $q_v$ denotes the
number of elements of the residue field cf. [Q2]. According to the
conjecture of Quillen and Lichtenbaum the map to the continuous
Galois cohomology $$K_{2n{+}1}(F)\otimes \Z_l\longrightarrow
H^1(G_{F,S_l},\, \Z_l(n{+}1))$$ constructed by Dwyer and
Friedlander in [DF], should be an isomorphism, for all odd $l.$
Here we have denoted by $G_{F,S_l}{=}G(F_{S_l}/F)$ the Galois
group of the maximal extension $F_{S_l}/F$ contained in ${\bar
F},$ which is unramified outside of the set $S_l$ of primes in
${\cal O}_F$ over $l.$
\medskip

In this paper we apply an axiomatic setup (see Section 2 for the
details) which embraces simultaneously the Mordell-Weil groups of
abelian varieties and K-groups of number fields.
Let ${\cal O}$ be a ring with unity.~We use the
general framework of the Mordell-Weil systems $\{B(L)\}_L$ of
finitely generated left ${\cal O}$-modules, indexed by finite
extensions $L/F$, which was developed in [BGK2] for the solution
of the support problem for {l}-adic representations.~One part of
the structure of the Mordell-Weil system is a collection of
well-behaved maps taking values in the Selmer groups defined by
Bloch and Kato [BK] in the Galois cohomology of $G({\bar F}/F).$
In the case of $E(F)$ these maps are induced by the classical
Kummer homomorphisms, whereas for K-groups they are given by the
Dwyer-Friedlander maps.
\medskip

Our main results concern linear dependence of nontorsion points in
$B(F)$. We assume that the system $\{B(L)\}_L$ and the Galois
representation in question meet the assumptions formulated in
Section 2.
\bigskip

\proclaim{Theorem A} [Theorem 2.9]

\noindent
Assume that the ring ${\cal O}$ is a finitely generated 
free $\Z$-module. Let $P$ and $P_1,\dots, P_r$ be nontorsion elements 
of $B(F)$ such that ${\cal O} P$ is a free ${\cal O}$ module and 
$P_1, P_2,\dots, P_r$ are linearly independent over ${\cal O}.$ Denote by 
$\Lambda$ the submodule of $B(F)$ generated 
by $P_1,P_2,\dots, P_r.$ Assume that $r_{v}(P) \in r_{v}(\Lambda)$ 
for almost all primes $v$ of $F.$ Then there exists a natural number 
$a$ such that $aP\in \Lambda.$
\endproclaim

In addition to the methods of [BGK2], we use in the proof of Theorem A 
the Kummer theory which was developed 
by Ribet in [Ri]. We prove that for 
Mordell-Weil systems, which come from Tate modules of some abelian varieties 
$A$ with $End (A) = \Z$ and from odd K-groups of number fields, 
a stronger result than Theorem A holds. Our Theorem 3.12 
shows that for such Mordell-Weil systems one can choose $a{=}1$  
in Theorem A. The main technical result of this part of the 
paper is Theorem 3.1. In the first step of the proof of Theorem 3.1 we 
applied an argument due to Khare and Prasad (cf. [KP], Lemma 5) 
extended to the context of the Mordell-Weil systems.
\medskip

The problem which we consider in this paper was motivated by the
support problem of Erd{\" o}s and also by our papers [BGK1] and
[BGK2] on the generalization of Erd{\" o}s problem to
${l}${-}adic Galois representations. The support problem for 
an abelian variety $A$ cf. [C-SR, p.277] is the following question: 
{\sl are the points $P,P_1\in A(F)$ related over the 
endomorphism algebra of $A,$ if the order of $r_v(P)\in A_v(\kappa_v)$ 
divides the order of $r_v(P_1),$ for almost all $v$ ?} Note that the groups 
$A_v(\kappa_v)$ are not necessarily cyclic, and therefore our 
problem for $A$ and $r{=}1,$ is 
essentially different from the support problem for $A.$ 
We would like to mention that Michael Larsen has recently given a solution 
of the support problem for all abelian varieties cf. [L].  
\medskip

In the case of abelian varieties our method provides the complete
solution of the problem of detecting linear dependence of
nontorsion points by reductions for the class of abelian varieties
for which Serre proved in [Se1] the analog of the open image
theorem.
\medskip

\proclaim{Theorem B} [Theorem 4.2]

\noindent
Let $A$ be a principally polarized abelian variety of dimension
$g$ defined over the number field $F$ such that $End(A) = \Z$
and $dim(A)=g$ is either odd or $g=2$ or $6.$
Let $P$ and $P_1,\dots, P_r $ be nontorsion elements of $A(F)$ such that $P_1,P_2,\dots, P_r$ are
linearly independent over ${\Z}.$ Denote by $\Lambda$ 
the subgroup of $A(F)$ generated 
by $P_1,P_2,\dots, P_r.$ Then the following two statements are equivalent:
\roster
\item[1]\,\,$P\in \Lambda$
\item[2]\,\,$r_{v}(P) \in r_{v}(\Lambda)$ for almost all primes $v$ of $F,$
where $r_v\colon A(F)\longrightarrow A_v(\kappa_v)$ are the reduction maps.
\endroster
\endproclaim

In particular, the Mordell-Weil group $E(F)$ of any non-CM
elliptic curve $E$ defined over the number field $F$ has the
property stated in Theorem C, i.e., our initial question has a
positive answer for such curves (cf. Corollary 4.3).~We hope that
in the case of elliptic curves defined over $\Q,$ Theorem B will
find some practical implementations.
\medskip

In the case of the K-group $K_{2n{+}1}(F),$ where $n\geq 0$ and $F$ is a 
number field, we put $B(F){=}K_{2n+1}(F)/C_{F},$ where $C_F$ is the 
subgroup of $K_{2n+1}(F)$ generated by $l$-parts, 
of kernels of the Dwyer-Friedlander maps, for all primes $l.$ 
Note that $C_F$ is a finite group by [DF], and if the Quillen-Lichtenbaum 
conjecture holds true, then $B(F){=}K_{2n{+}1}(F)$ up to $2$-torsion.  
\bigskip

\proclaim{Theorem C}[Theorem 4.1]

\noindent 
Let $P$ and $P_1,$ $P_2,$ $\dots,$ $P_r$
be nontorsion elements of $K_{2n+1}(F)/C_{F},$ such that 
$P_1,$ $P_2,$ $\dots,$ $P_r$ are
linearly independent over ${\Z}.$ Let $\Lambda\subset K_{2n+1}(F)/C_{F}$ 
denote the subgroup generated by $P_1,P_2,\dots, P_r.$ The following two 
statements are equivalent: 
\roster 
\item[1]\,\,$P\in \Lambda$ 
\item[2]\,\,$r_{v}(P) \in r_{v}(\Lambda)$ for almost all primes $v$ of $F.$
\endroster
\endproclaim

For two points, i.e., if $r{=}1,$ Theorem C was proven already in [BGK1].
If $n=0,$ then Theorem C gives the statement about the multiplicative 
group of the field $F,$ which was previously proven by Khare, [K], 
Proposition 3, p.10.
\medskip

While this paper was being prepared we have learned that Tom
Weston has proven by a different
method a result similar to our Theorem 2.9, for all 
abelian varieties with commutative algebras of endomorphisms 
cf. [We], Corollary. 2.8. 
\bigskip\bigskip

\subhead 2. Kummer theory for $l$-adic representations
\endsubhead

\subsubhead{Notation}
\endsubsubhead

\roster
\item"{$l$}"\quad is a prime number
\item"{$F$}" \quad is a number field, ${\cal O}_{F}$ its ring of integers
\item"{$F_S$}" \quad is the maximal extension of $F$ unramified outside
a finite set $S$ of primes 
\item"{}"\quad in ${\cal O}_{F}$
\item"{$G_F$}"\quad $=\quad G({\bar F}/F)$
\item"{$G_{F,S}$}"\quad $=\quad G(F_{S}/F)$
\item"{$v$}"\quad denotes a finite prime of ${\cal O}_{F},$ $\kappa_v{=}{\cal O}_F/v$ denotes the 
residue field at $v$
\item"{$g_v$}" \quad $=\quad G({\overline \kappa_v}/\kappa_v)$
\item"{$T_l$}"\quad denotes a free ${\Z_l}$-module of finite rank $d$
\item"{$V_l$}" \quad $=\quad T_l\otimes_{\Z_l} \Q_l$
\item"{$A_l$}" \quad $= \quad V_l/T_l$
\item"{$\rho_l:$}"\quad $G_F\rightarrow GL(T_l)$ is a Galois representation
unramified outside a fixed finite
\item"{}"\quad set $S_l$ of primes of ${\cal O}_{F}$ containing all primes
above $l$
\item"{${\overline{\rho_{l^k}}}$}" \quad denotes the residual representation
 $G_F \rightarrow GL(T_l/{l^k})$ induced by $\rho_l$
\item"{$F_{l^k}$}"\quad $=\quad F(A_l[l^k]),$ for any $k>0,$ denotes the number field
${\bar F}^{ker \overline{\rho_{l^k}}}$
\item"{$F_{l^\infty}$}"\quad $=\quad {\bigcup}_{k}F_{l^k}$
\item"{$G_{l^k}$}"\quad $=\quad G(F_{l^k}/F)$
\item"{$G_{l^{\infty}}$}"\quad $=\quad G(F_{l^{\infty}}/F)$
\item"{$H_{l^k}$}"\quad $=\quad G({\overline F}/F_{l^k})$
\item"{$H_{l^{\infty}}$}"\quad $=\quad G({\overline F}/F_{l^{\infty}})$
\item"{$C[l^k]$}"\quad denotes the subgroup of $l^k${-}torsion elements of an abelian group $C$
\item"{$C_l$}"\quad  $=\quad \bigcup_{k\geq 1}C[l^k],$ is the $l${-}torsion subgroup of $C.$
\endroster

\noindent
Let $L/F$ be a finite extension and $w$ a finite prime in $L.$
To indicate that $w$ is not over any prime in $S_l$ we will write
$w\notin S_l,$ slightly abusing notation.
Let $\cal O$ be a ring with unity which acts on $T_{l}$ 
in such a way that the action commutes 
with the $G_F$-action. Unless specified otherwise (see paragraph after 
Lemma 3.11), all modules over the ring 
${\cal O},$ which we consider in this paper, are left ${\cal O}$-modules. Let 
$\{B(L)\}_L$ be a direct system of finitely generated ${\cal O}$-modules
indexed by all finite field extensions $L/F.$ The structure maps of 
the system are induced by
inclusions of fields. We assume that for every embedding
$L \rightarrow L^{\prime}$ of extensions of $F,$ the structure map
$B(L) \rightarrow
B(L^{\prime})$ is a homomorphism of ${\cal O}$-modules.
Similarly, for every prime $v$ of $F$ we define a direct system
$\{B_{v}(\kappa_w)\}_{\kappa_{w}}$ of finite ${\cal O}$-modules 
where $\kappa_w$ is a residue
field for a prime $w$ over $v$ in a finite extension $L/F.$
We require that $G_{F}$ acts on both systems:  $\{B(L)\}_L$ and
$\{B_{v}(\kappa_w)\}_{\kappa_{w}}$ functorially. Let us put $B({\overline
F}) = \varinjlim_{L/F}\, B(L).$ We assume that the actions of $G_{F}$ and
${\cal O}$ have the following properties:
\medskip

\roster
\item"{($A_{1}$)}"  for each $l,$ each
finite extension $L/F$ and any prime $w$ of $L,$ such
that $w\not\in S_l,$ we have $T_l^{Fr_w} = 0,$ where $Fr_w\in g_w$ denotes
the arithmetic Frobenius at $w.$
\endroster
\roster
\item"{($A_{2}$)}"  for every $L$ and $w\notin S_{l}$ there are natural maps
${\psi_{l,L}},$ ${\psi_{l,w}}$ and $r_w$
respecting $G_{F}$ and $\cal O$ actions such that the diagram
commutes:
$$
\CD
B(L)\otimes \Z_l @>{r_w}>> B_{v}(\kappa_w)_{l}\\
@V{\psi_{l,L}}VV @V{\cong}V{\psi_{l,w}}V\\
H^{1}_{f,S_{l}}(G_{L},T_{l}) @>{r_{w}}>> H^{1}(g_{w},T_{l})
\endCD
\tag{2.1}$$ where $H^{1}_{f,S_{l}}(G_{L},T_{l})$ is the group
defined by Bloch and Kato (cf.[BK], see also [BGK2]). The left
(resp., the right) vertical arrow in the diagram (2.1) is an
embedding (resp., an isomorphism) for every $L$ (resp., for every
$w\notin S_{l}$)
\endroster
\roster
\item"{($A_{3}$)}" for every $L$ the map ${\psi_{l,L}}$ is an isomorphism
for almost all $l$ or $B({\overline F})$ is a discrete $G_F$-module
divisible by $l,$ for almost all $l,$ and for every $L$ we have:
$B({\overline F})^{G_L} \cong B(L)$ and $H^0(G_L; A_l) \subset B(L).$
\endroster

\noindent For a point $R\in B(L)$ (resp., a subgroup $\Gamma
\subset B(L)$) we denote $\hat R = \psi_{l,L}(R)$ (resp.
$\hat\Gamma {=}\psi_{l,L}(\Gamma)$). As in [Ri] we impose the
following four axioms on the representations which we consider:
\medskip

\roster
\item"{($B_1$)}"\quad  $End_{G_{l}} (A_{l}[l])\cong {\cal O}/l{\cal O},$ for almost all $l$ and
\item"{}"\quad $End_{G_{l^{\infty}}}(T_l) \cong {\cal O}\otimes {\Z_l},$ for all $l$
\endroster
\roster
 \item"{($B_2$)}"\quad  $A_{l}[l]$ is a semisimple
$\F_{l}[G_{l}]${-}module, for almost all $l$ and
\item"{}"\quad $V_{l}$ is a semisimple
${\Q}_{l}[G_{l^{\infty}}]${-}module, for all $l$
\endroster
\roster
\item"{($B_3$)}"\quad  $H^{1}(G_{l}; A_{l}[l]) = 0,$ for almost all $l$ and
\item"{}"\quad $H^{i}(G_{l^{\infty}}; T_{l})$ are finite groups, for all $l$ and all $i \geq 0$
\endroster
\roster
\item"{($B_4$)}"\quad for each finitely generated subgroup $\Gamma \subset B(F)$ the group
$$\Gamma^{\prime} = \{ P\in B(F) : mP \in \Gamma \quad \text{ for
some }\quad m\in \N \}$$
\quad is such that $\Gamma^{\prime}/{\Gamma}$ has a finite exponent.
\endroster

\noindent
Note that in the case of the Tate module of an abelian variety,
the conditions ($A_1$)-($A_3$) were checked in [BGK2], Section 3.
The conditions ($B_1$)-($B_4$) are fulfilled in this case due to results 
of Faltings, Zarhin, Serre, Mordell and Weil (cf. the proof of Theorem 4.2). In the sequel we assume that ${\cal O}$ is a finitely generated free
$\Z$-module.
\medskip

\noindent
The next lemma describes the relation between the conditions in the 
assumption ($B_1$).

\proclaim{Lemma 2.2} If $End_{G_{l}} (A_{l}[l]) \cong {\cal
O}/l{\cal O},$ then $End_{G_{l^{\infty}}}(T_l)\cong {\cal
O}\otimes {\Z_l}.$
\endproclaim

\demo{Proof} We prove by induction that for every $k>0$ we have
$End_{G_{l^{k}}}(A_{l}[l^k]) \cong {\cal O}/{l^k}{\cal O}.$ Assume
that the claim is true for $k{=}1, 2, \dots, i{-}1.$ Consider the
commutative diagram: $$ \eightpoint \CD 0 @>>>
End_{G_{l^{i{-}1}}}(A_{l}[l^{i{-}1}]) @>>>
End_{G_{l^{i}}}(A_{l}[l^{i}]) @>{p}>> End_{G_{l}}(A_{l}[l])@>>>
0\\ @. @A{\cong}AA @A{s}AA @A{\cong}AA @.
\\ 0 @>>> {\cal O}/l^{i{-}1} {\cal O}  @>>> {\cal O}/l^{i} {\cal O}
@>{p^{\prime}}>> {\cal O}/l {\cal O}
@>>> 0.
\endCD$$
The map $p$ in the diagram is a surjection, since $p^{\prime}$ is a surjection. Hence, $s$ is
an isomorphism. The lemma follows by taking inverse limit over $k.$ \qed
\enddemo
\medskip

Let $\Lambda$ be a finitely generated $\cal O${-}submodule of $B(F).$ Throughout the paper we
assume that the points $P_1,\dots ,P_r$ constitute a basis of $\Lambda $ over ${\cal O}$ i.e., they
give an ${\cal O}${-}isomorphism $\Lambda \cong {\Cal O}^r.$ Let $\overline P_i$ denote
 the image of $P_i$ in $\Lambda/l\Lambda.$ It is clear that ${\overline P_1},\dots ,\overline
 P_r$ form a basis of the module $\Lambda/l\Lambda$ over ${\cal O}/l{\cal O}.$
\medskip

In the remainder of this section, following [Ri] we introduce the
Kummer theory for the $l${-}adic representations which meet our
assumptions. For $P \in B(F)$ and $k>0$ we have the Kummer maps:
$${{\phi}^{(k)}_{P}} : H_{l^{k}} \longrightarrow
A_{l}[l^{k}]\tag{2.3}$$ $${{\phi}^{(k)}_{P}}(\sigma)= \sigma
(\frac{1}{l^k} {\hat P})-{\frac{1}{l^k}}{\hat P}.$$ We define:
$$\Phi^{(k)}:H_{l^{k}} \longrightarrow
\bigoplus_{i=1}^{r}A_{l}[l^k]$$ $${\Phi}^{(k)}
=({\phi}^{(k)}_{P_1},...,{\phi}^{(k)}_{P_r}).$$ One checks easily
that, for every $k>1,$ the following diagram commutes: $$ \CD
H_{l^{k}} @>{\phi}^{(k)}_{P}>>  A_{l}[l^{k}]\\ @VVV @VV{\times
l}V\\H_{l^{k-1}}  @>{\phi}^{(k-1)}_{P}>>
 A_{l}[l^{k{-}1}].
\endCD
\tag{2.4}$$ Taking the inverse limits over $k$ in (2.4) we obtain
a map: $${{\phi}^{(\infty)}_{P}} : H_{l^{\infty}} \longrightarrow
T_{l}(A),$$ which we denote briefly
${\phi}_{i}:={\phi}^{(\infty)}_{P_{i}}.$ We define:
$$\Phi:H_{l^{\infty}} \longrightarrow \bigoplus_{i=1}^{r}T_l$$
$${\Phi}=({\phi}_{1},\dots,{\phi}_{r}).$$
\medskip

\noindent
Consider the field $F_{l}(\frac{1}{l}\hat\Lambda):={\overline F}^{ker\, \Phi^{(1)}}$
and the map induced by $\Phi^{(1)}$:
 $$G(F_{l}(\frac{1}{l}\hat\Lambda)/F_{l}) \longrightarrow \bigoplus_{i=1}^{r} A_{l}[l]\tag{2.5}$$
cf. [BGK2], (4.13). Since $\overline P_1,\dots, \overline P_r$ is a basis of
$\Lambda/l\Lambda$ over ${\cal O}/l{\cal O},$ the same argument
as in [Ri],  Theorem 1.2, shows that, for every $l,$ the map (2.5) is
an isomorphism of $G_{l}${-}modules.

\proclaim{Lemma 2.6} For almost all $l,$ there exists an
isomorphism of ${\cal O}/l{\cal O}${-}modules: $$\Lambda /l\Lambda
\longrightarrow Hom_{G_l}(G(F_{l}(\frac{1}{l}\hat\Lambda)/F_{l});
A_{l}[l]). \tag{2.7}$$
\endproclaim
\medskip

\demo{Proof} By Kummer theory the map $P\mapsto \phi_P$ induces a
homomorphism: $$ B(F)/lB(F)\longrightarrow Hom_{G_l}(H_l;\,
A_{l}[l])=H^1(H_l;\, A_{l}[l])^{G_l} $$ which fits into the
commutative diagram: $$ \diagram {\phantom{\Big|}
\Lambda/l\Lambda\phantom{\Big|} }  \ar@{^{(}->}[d] \ar[r]&
{\phantom{\Big|} Hom_{G_l}(G(F_{l}(\frac{1}{l}\hat\Lambda)/F_{l});
A_{l}[l])\phantom{\Big|} } \ar@{^{(}->}[d]\\ {\phantom{\Big|}
B(F)/lB(F)\phantom{\Big|}} \ar[r]\ar@{^{(}->}[d] &
{\phantom{\Big|} Hom_{G_l}(H_l;\, A_{l}[l])\phantom{\Big|}
}\ar_{=}[d]\\ {\phantom{\Big|} H^1(G_F;\, A_{l}[l])\phantom{\Big|}
} \ar[r]^{res}& {\phantom{\Big|} H^1(H_l;\,
A_{l}[l])^{G_l}\phantom{\Big|} } \\
\enddiagram
 \tag{2.8}$$
as the middle horizontal map. The map (2.7) is the upper
horizontal map in the diagram (2.8). The lower horizontal map is
the restriction. By the inflation-restriction sequence in Galois
cohomology, the kernel of $res$ is contained in $H^1(G_l;\,
A_{l}[l])$ which vanishes, for almost all $l,$ by axiom
($B_3$). The left-lower vertical map in (2.8) is injective by
($A_3$). The left-upper vertical arrow is an injection for almost
all $l$ by ($B_4$). This shows that the map (2.7) is an imbedding
for almost all $l.$ Comparing the dimensions over $\Z /l{\Z},$ we
conclude by axiom ($B_1$) and the isomorphism (2.5), that the
map (2.7) is an isomorphism of ${\cal O}/l{\cal O}${-}modules, for
almost all $l.$ \qed\enddemo
\medskip

\proclaim{Theorem 2.9}

\noindent
Let $P\in B(F)$ be such that ${\cal O} P$ is a free ${\cal O}$ module and 
for almost all primes $v $ of $F,$ we have $r_{v}(P) \in r_{v}(\Lambda).$ Then
there is a natural number $a$ such that $aP\in \Lambda.$
\endproclaim
\medskip

\demo{Proof}
Let $\Lambda_{1}$ be the $\cal O${-}submodule of $B(F)$ generated by $P$. We consider the Galois extensions
$F_{l}(\frac{1}{l}\hat\Lambda)/F$ and $F_{l}(\frac{1}{l}\hat\Lambda_{1})/F,$ where
$F_{l}(\frac{1}{l}\hat\Lambda_{1}):= {\overline F}^{ker\, \phi_P^{(1)}}.$
Assume that $v$ splits completely in $F_{l}(\frac{1}{l}\hat\Lambda)/F.$ Let $w^{\prime}$ be any prime over $v$ in
 $F_{l}(\frac{1}{l}\hat\Lambda)/F.$ We have $\kappa_{w^{\prime}}=\kappa_v,$ hence $g_{w^{\prime}}=g_v.$ Therefore
$r_{w^{\prime}}(\frac{1}{l}\hat\Lambda)\in H^{1}(g_{v},T_{l}).$ Let
$w$ be any prime of $F_{l}(\frac{1}{l}\hat\Lambda_{1})$ over $v.$
 Hence, by the assumption $r_{w}(\frac{1}{l}\hat\Lambda_{1})\in H^{1}(g_{v},T_{l}).$  Let ${\hat R}$ be any point in
$\frac{1}{l}\hat{\Lambda}_{1}.$ Since $Fr_{w}(\hat R) = \hat R +
{\hat P_{0}},$ for some ${\hat P_{0}} \in A_{l}[l],$ we have: $$
r_{w}(\hat R)= Fr_{w}  r_{w}(\hat R)= r_{w}Fr_{w}(\hat
R)=r_{w}(\hat R)+r_{w} ({\hat P}_{0}).$$
By [BGK2], Lemma 2.13 and the axiom ($A_1$) we obtain $\hat P_{0}=0.$
This means that $v$ splits completely in the field $F_{l}(\frac{1}{l}\hat\Lambda_{1}).$
 By the Frobenius density theorem [J], Corollary 5.5, p.136, we get:
$$F_{l}(\frac{1}{l}\hat\Lambda_{1})\subset
F_{l}(\frac{1}{l}\hat\Lambda).\tag{2.10}$$ Hence, by the
isomorphism (2.7) applied to $\Lambda$ and $\Lambda_1,$ and the
inclusion (2.10) we obtain the diagram: $$ \CD \Lambda_{1}
/l\Lambda_{1} @>>> Hom_{G_l}(
G(F_{l}(\frac{1}{l}\hat\Lambda_{1})/F_{l}); A_{l}[l])\\ @VVV
@VV{\alpha}V\\ \Lambda /l\Lambda @>>> Hom_{G_l}(
G(F_{l}(\frac{1}{l}\hat\Lambda)/F_{l}); A_{l}[l]),
\endCD
\tag{2.11}$$

\noindent
where the left vertical arrow is defined to make the diagram (2.11) commute. Since $\alpha$ is injective by (2.10),
we obtain an inclusion:
$$\Lambda_{1} /l\Lambda_{1} \subset \Lambda /l\Lambda.$$
Put $B_{0} = B(F)/{\Lambda}.$ By (2.11) we see that $P$ maps to zero in $B_{0}/l{B_{0}}.$
Hence, by the axiom ($B_4$) we get that $P{-}Q \in (B_{0})_{tor},$ for some $Q\in\Lambda .$
This proves the theorem, if we take for $a$ the exponent of the finite group $(B_{0})_{tor}.$\qed\enddemo

\proclaim{Lemma 2.12} If ${\alpha}_{1},\dots ,{\alpha}_{r}\in {\Cal
O}\otimes_{\Z}{{\Z}}_{l}$ are such that
${{\alpha}_{1}}{\phi_{1}}+\dots +{{\alpha}_{r}}{\phi_{r}}=0,$ then
${{\alpha}_{1}}=\dots ={{\alpha}_{r}}=0.$
\endproclaim

\demo{Proof} Let $\Psi$ be the composition of maps:
$$B(F)\otimes_{\Z}{\Z}_{l}\hookrightarrow
H^{1}(G_{F};T_{l})\longrightarrow H^{1}(H_{l^{\infty}};T_{l}).$$
Note that $H^{1}(H_{l^{\infty}};T_{l})=Hom(H_{l^{\infty}};T_{l})$
and $\Psi(P_{i}\otimes 1)=\phi_{i}.$ By ($B_{3}$) we have
$ker\Psi\subset (B(F)\otimes_{\Z}{\Z}_{l})_{tor}.$ Let
$t:=\#B(F)_{tor}.$ Since $\Psi$ is an ${\Cal
O}{\otimes}_{\Z}{{\Z}}_{l}${-}homomorphism, we have: $$0\quad
=\quad {{\alpha}_{1}}{\phi_{1}}+\dots
+{{\alpha}_{r}}{\phi_{r}}\quad =\quad
\Psi({{\alpha}_{1}}(P_{1}\otimes 1)+\dots
+{{\alpha}_{r}}(P_{r}\otimes 1)).$$ Hence,
${{\Sigma}_{j=1}^{r}}{{\alpha}_{j}}(P_j\otimes 1)\in
(B(F)\otimes_{\Z}{\Z}_{l})_{tor},$ so:
$$t{{\alpha}_{1}}({P_{1}\otimes 1)}+\dots
+t{{\alpha}_{r}}({P_{r}\otimes 1}) =0$$ in
$B(F)\otimes_{\Z}{\Z}_{l}.$ Observe that the points
$P_1{\otimes}1,\dots ,P_r{\otimes}1$ are linearly independent over
$\cal O \otimes_{\Z}{\Z}_{l}$ in $B(F)\otimes_{\Z}{\Z}_{l}.$ This
implies that $t{{\alpha}_{j}}=0$ for $ j=1,\dots, r.$  Hence,
$${{\alpha}_{1}}=\dots ={{\alpha}_{r}}=0,$$ because ${\Cal O}$ is
a free ${\Z}${-}module, by  assumption. \qed\enddemo

\proclaim{Lemma 2.13}The image of the map $\Phi$ is
 open in $\bigoplus_{i=1}^{r}T_l$ in the $l$-adic topology.

\endproclaim

\demo{Proof} It is enough to show that $Im\,\Phi$ has a finite index in
$\bigoplus_{i=1}^{r}T_l.$ The proof follows the lines of the proof of [Ri, Th. 1.2].
Let $W=\bigoplus_{i=1}^{r}V_l$ and $M=Im(\Phi\otimes 1)
\subset W,$ where $V_l{=}T_l{\otimes}_{\Z_{l}}{\Q_{l}}$. Both $M$
and $W$ are ${\Q}_{l}[G_{l^{\infty}}]${-}modules.
First we will show that $\Phi\otimes 1$ is onto. Suppose it is not.
By $(B_{2})$ we have a nontrivial decomposition $W=M\oplus M_{1}$ of ${\Q}_{l}[G_{l^{\infty}}]${-}modules.
Let $\pi_{i}:W\rightarrow V_{l}$ be a projection that
maps $M_{1}$ nontrivially. By the axiom ($B_1$) we have:
$$\pi_{i}(v_{1},\dots ,v_{r})=\sum_{j=1}^{r} {{\beta}_{j}}v_{j},$$
for some ${{\beta}_{j}}\in {\Cal O}\otimes \Q_l.$ Since $\pi_{i}$ is nontrivial, we see that some ${{\beta}_{j}}$ is
nonzero. On the other hand
$${\pi_{i}}({\Phi}\otimes 1)(h)=\sum_{j=1}^{r} {{\beta}_{j}}({\phi}_{j}(h)\otimes 1)=0,$$
for all $h\in H_{l^{\infty}}.$ Since ${{\beta}_{j}}\in {\Cal O}\otimes \Q_l,$ we can multiply the
last equality by a suitable power $l^k$ to get:
$$0=\sum_{j=1}^{r}{{\alpha}_{j}}({\phi}_{j}(h)\otimes 1),$$
where
${{\alpha}_{j}}=l^k \beta_j \in {\Cal O}\otimes \Z_l.$
Since the maps:
$${\Cal O}\otimes \Z_l\hookrightarrow {\Cal O}\otimes \Q_l$$
$$Hom(G_{l^{\infty}},T_l)\hookrightarrow Hom(G_{l^{\infty}},V_l)$$
are imbeddings of ${\Cal O}\otimes\Z_l${-}modules, we obtain: $\sum_{j=1}^{r}{{\alpha}_{j}}{\phi}_{j}=0.$
By Lemma 2.12 this shows:
$${{\alpha}_{1}}=\dots ={{\alpha}_{r}}=0,$$
 a contradiction. Therefore we have $M_1=0.$ Since both $\bigoplus_{i=1}^{r}T_l$ and $Im\,\Phi$
are $\Z_l${-}lattices in $\bigoplus_{i=1}^{r}V_l,$ we see that $Im\,\Phi$ has a finite index in
$\bigoplus_{i=1}^{r}T_l,$ as required. \qed\enddemo
\bigskip

\subhead{3. Main Results}
\endsubhead

In this section we show that for some $\{B(L)\}_L,$ one can choose
$a{=}1$ in Theorem 2.9. Let ${\cal G}_{l}^{alg}$ be the Zariski
closure of the image of $\rho_l$ in the algebraic group scheme
\newline
$GL_d/{\Z}_l,$ endowed with the unique structure of the reduced,
closed group subscheme. The next proposition is the main technical
result of the paper.
\medskip

\proclaim{Theorem 3.1}

\noindent
 Assume that $\rho_{l}(G_{F})$ is open
in ${\cal G}_{l}^{alg}(\Z_{l})$ in the topology induced from $\Z_l.$ In addition, assume
that $\rho_{l}(G_{F})$ contains an open subgroup of the
group of homotheties and that the reduction map ${\cal
G}_{l}^{alg}(\Z_{l})\rightarrow {\cal G}_{l}^{alg}(\Z/l^{k})$ is
onto, for every $k>0.$ Let $P_1,\dots, P_r\in B(F)$ be points of infinite order, which are linearly independent
over ${\cal O}.$ Let
$$I = \{i_1,\dots, i_s\} \subset \{1, \dots, r\}$$
 be any subset of indices and let
$$J = \{j_1, \dots, j_{r-s}\} \subset \{1, \dots, r\}$$
be such that $I \cap J = \emptyset$ and $I \cup J = \{1, \dots, r\}.$ Then for any natural $M$, and for any prime $l,$
there are infinitely many primes $v,$ such that the images of the points
$P_{i_1},\dots, P_{i_s}$ via the map:
$$r_{v}\,:\,B(F)\rightarrow B_v(\kappa_v)_{l}$$
are trivial and the images of the points $P_{j_1},\dots, P_{j_{r{-}s}}$
have orders divisible by $l^M.$
\endproclaim

\demo{Proof} {\bf Step 1.} This part  of the proof is analogous to
the proof of Lemma 5 of [KP]. Let ${\Lambda}_{I}$ (resp.,
$\Lambda_J$) be the ${\cal O}${-}submodule of ${\Lambda}$
generated by $P_{i_1},\dots,P_{i_r}$ (resp., by $P_{j_1},\dots,
P_{j_{r{-}s}}$). Let
$F_{l^{\infty}}(\frac{1}{l^{\infty}}\hat{\Lambda}_{I})$ (resp.,
$F_{l^{\infty}}(\frac{1}{l^{\infty}}\hat{\Lambda}_{J})$) be the
composite of the fields
$F_{l^k}(\frac{1}{l^k}{\hat{\Lambda}_{I}})$ (resp.,
$F_{l^k}(\frac{1}{l^k}{\hat{\Lambda}_{J}})$), for $k>0.$ We will
use the notation similar to [KP]:
\medskip

\roster
\item"{${\cal D}_k$}"\quad $=$\quad $G(F_{l^k}(\frac{1}{l^k}{\hat{\Lambda}_{I}})/{F}_{l^{k}})$
\item"{${\cal E}_k$}"\quad $=$\quad $ G(F_{l^k}( \frac{1}{l^k}\hat{\Lambda}_{I})/F)$
\item"{${\cal D}_{\infty}$}"\quad $=$\quad $G(F_{l^{\infty}}(\frac{1}{l^{\infty}}{\hat{\Lambda}_{I}})/
{F}_{l^{\infty}})$
\item"{${\cal E}_{\infty}$}"\quad $=$ \quad $G(F_{l^{\infty}}(\frac{1}{l^{\infty}}\hat{\Lambda}_{I} )/F).$
\endroster
\medskip

\noindent
Consider the following commutative diagram:
$$\diagram
 0 \ar[r]& {\phantom{\Big|} {\cal D}_{{\infty}}
\phantom{\Big|}} \ar[r]\ar@{^{(}->}[d]&{\phantom{\Big|}
{\Cal E}_{{\infty}}\phantom{\Big|} }\ar[r]\ar@{^{(}->}[d]
& \ar[r]\ar@{^{(}->}[d] {\phantom{\Big|} {G}_{l^{\infty}}
\phantom{\Big|} }& 0\\
0 \ar[r] & {\phantom{\Big|}
({{\G}^{d}_{a}})^{r}(\Z_{l})\phantom{\Big|}} \ar[r]&
{\phantom{\Big|} {\Cal E}(\Z_{l})\phantom{\Big|}}
\ar[r]& {\phantom{\Big|} {\Cal G_{l}^{alg}}(\Z_{l})\phantom{\Big|}}
\ar[r]& 0.
\enddiagram
\tag{3.2}$$

\noindent
The group scheme ${\Cal E}= {\Cal G}^{alg}_{l} \ltimes ({{\G}^{d}_{a}})^{r}$
is a semi-direct product and $({{\G}^{d}_{a}})^{r}(\Z_{l})\cong T_{l}^{r}:=\bigoplus_{i{=}1}^rT_l.$
The left vertical arrow in the diagram (3.2) is induced by the
 map $\Phi.$ The image of $\Phi$ is open by Lemma 2.13. Since by 
assumption $G_{l^{\infty}}$ is open in ${\Cal G_{l}^{alg}}(\Z_{l}),$ it follows by (3.2) that ${\Cal E}_{{\infty}}$
is open in ${\Cal E}(\Z_{l}).$ Because the map ${\cal G}_{l}^{alg}(\Z_{l})\rightarrow {\cal G}_{l}^{alg}(\Z/l^{k})$ is
onto, the map ${\cal E}(\Z_{l})\rightarrow {\cal E}(\Z/l^{k})$ is also onto, for $k>0.$
Consider the following commutative diagram.

$$\eightpoint \diagram
&{\phantom{\Big|} (\Z/l^{k+1})^{dr}\phantom{\Big|}}\ar@{->>}'[d][dd]\ar@{^{(}->}[rr]&&{\phantom{\Big|}
{\cal E}(\Z /l^{k+1})\phantom{\Big|}}\ar@{->>}'[d][dd]\ar@{->>}[rr]
&&{\phantom{\Big|} {\cal G}_{l}^{alg}(\Z/l^{k+1})\phantom{\Big|}}\ar@{->>}[dd]\\
 {\phantom{\Big|} {\cal D}_{{k{+}1}}\phantom{\Big|}}\ar[dd]\ar@{^{(}->}[rr]\ar@{^{(}->}[ur]
&&{\phantom{\Big|} {\cal E}_{{k{+}1}}\phantom{\Big|}}\ar@{->>}[dd]\ar@{->>}[rr]\ar@{^{(}->}[ur]&&
{\phantom{\Big|} {G}_{l^{k{+}1}}\phantom{\Big|}}\ar@{->>}[dd]\ar@{^{(}->}[ur]&&\\
&{\phantom{\Big|} (\Z/l^{k})^{dr}\phantom{\Big|}}
\ar@{^{(}->}'[r][rr]&&{\phantom{\Big|} {\cal E}(\Z/l^{k})\phantom{\Big|}}\ar@{->>}'[r][rr] &&
{\phantom{\Big|} {\cal G}_{l}^{alg}(\Z/l^{k}) \phantom{\Big|}}\\
{\phantom{\Big|}{\cal D}_{{k}}\phantom{\Big|} }\ar@{^{(}->}[rr]\ar@{^{(}->}[ur] &&
{\phantom{\Big|} {\cal E}_{{k}}\phantom{\Big|} }\ar @{->>}[rr]\ar@{^{(}->}[ur]&&
{\phantom{\Big|} {G}_{l^{k}}\phantom{\Big|}}\ar@{^{(}->}[ur]&&
\enddiagram\tag{3.3}$$

\noindent Using the argument on the natural congruence subgroup of
${\cal E}(\Z_{l})$ of level $l^k,$ as in  the proof of Lemma 5,
[KP],  we note that, for $k$ big enough, the preimage of ${\cal
E}_{{k}}$ via the middle vertical arrow in the rear wall of
diagram (3.3) is ${\cal E}_{{k{+}1}}.$ It follows by diagram
chasing that the left vertical arrow ${\cal D}_{{k{+}1}}
\rightarrow {\cal D}_{{k}}$ in the front wall of (3.3) is
surjective. This immediately implies that, for $k$ big enough:
$$F_{l^k}(\frac{1}{l^k}{\hat{\Lambda}_{I}})\,\,\cap\,\,
F_{l^{k+1}}\quad =\quad F_{l^k}.\tag{3.4}$$
\medskip

\noindent {\bf Step 2.}  We will make use of the following tower
of fields. $$\diagram
&&F_{l^{k+1}}(\frac{1}{l^{k}}{\hat{\Lambda}_{I}}
,\frac{1}{l^{k}}{\hat{\Lambda}_{J}})&\\
&F_{l^{k}}(\frac{1}{l^{k}}{\hat{\Lambda}_{I}},\frac{1}{l^{k}}{\hat{\Lambda}_{J}})\urline
&F_{l^{k+1}}(\frac{1}{l^{k}}{\hat{\Lambda}_{I}})\uline&&
\\ &F_{l^{k}}(\frac{1}{l^{k}}{\hat {\Lambda}_{I}})\urline\uline^{\sigma}
&& F_{l^{k+1}}\ulline&&\\ &&F_{l^{k}}\ulline_{id}\urline^{h}&&\\
&&F\uline&&&&\\
\enddiagram\tag{3.5}
$$

\noindent
Consider the following commutative diagram:
$$ \diagram 0  \ar[r]
&{\phantom{\Big|} G_1\phantom{\Big|}}  \ar[r]{\ar@{^{(}->}[d]}& {\phantom{\Big|} G_2\phantom{\Big|}} \ar[r]\ar@{^{(}->}[d]& {\phantom{\Big|} G_3\phantom{\Big|}}
\ar[r]\ar@{^{(}->}[d]\ar[r]& 0\\
 0 \ar[r]& {\phantom{\Big|} T_l^{r{-}s}\phantom{\Big|}} \ar[r] & {\phantom{\Big|} T_l^r\phantom{\Big|}}
\ar[r]& {\phantom{\Big|} T_l^{s}\phantom{\Big|}} \ar[r] &  0,\\
\enddiagram
 \tag{3.6}$$
where we have put:

\roster
\item"{$G_1$}"\quad  $=$ $G(F_{l^{\infty}}(\frac{1}{l^{\infty}}\hat\Lambda _{I},\frac{1}{l^{\infty}}\hat\Lambda_{J})
/F_{l^{\infty}}(\frac{1}{l^{\infty}}\hat\Lambda_{I}))$
\item"{$G_2$}"\quad  $=$ $G(F_{l^{\infty}}(\frac{1}{l^{\infty}}\hat\Lambda_{I},
\frac{1}{l^{\infty}}\hat\Lambda_{J})/F_{l^{\infty}})$
\item"{$G_3$}"\quad  $=$ $G(F_{l^{\infty}}(\frac{1}{l^{\infty}}\hat\Lambda
_I)/F_{l^{\infty}})$
\endroster

\noindent and $F_{l^{\infty}}(\frac{1}{l^{\infty}}\hat\Lambda
_{I},\frac{1}{l^{\infty}}\hat\Lambda_{J})=
F_{l^{\infty}}(\frac{1}{l^{\infty}}\hat\Lambda_{I})F_{l^{\infty}}(\frac{1}{l^{\infty}}\hat\Lambda_{J}).$
All vertical arrows in the diagram (3.6) are the Kummer maps
discussed in Section 2.~The right and the middle vertical arrows
have open images (equivalently, finite cokernels) by Lemma 2.13.
Applying the snake Lemma to the diagram (3.6), we observe that the
left vertical arrow has finite cokernel, hence it has an open
image. Consider the following commutative diagram:

$$ \diagram & G(F_{l^{\infty}}(\frac{1}{l^{\infty}}\hat\Lambda_I,
\frac{1}{l^{\infty}}\hat\Lambda_J)/F_{l^{\infty}}(\frac{1}{l^{\infty}}\hat\Lambda_I))
 \quad\ar@{^{(}->}[r]\ar[d]&\quad T_l^{r{-}s} \ar@{->>}[d]&\\
 & G(F_{l^{k}}(\frac{1}{l^{k}}\hat\Lambda_I,
\frac{1}{l^{k}}\hat\Lambda_J)/F_{l^{k}}(\frac{1}{l^{k}}\hat\Lambda_I))\quad
\ar@{^{(}->}[r] &\quad \bigoplus_{i{=}1}^{r{-}s}A_l[l^k].&&\\
\enddiagram \tag{3.7}$$

\noindent
The horizontal arrows in the diagram (3.7) are the Kummer maps.
The upper horizontal arrow is the left vertical arrow in the diagram (3.6), so it has
an open image. Hence, the lower horizontal arrow in the diagram (3.7) has the cokernel
bounded independently of $k.$ Let $M$ be a natural number. It follows that,
for $k$ big enough, there is an element $\sigma \in G(F_{l^{k}}(\frac{1}{l^{k}}\hat\Lambda_I,
\frac{1}{l^{k}}\hat\Lambda_J)/F_{l^{k}}(\frac{1}{l^{k}}\hat\Lambda_I)),$
such that $\sigma$ maps via the horizontal arrow in (3.7)
to an element of $\bigoplus_{i{=}1}^{r{-}s} A_l[l^k],$ with all 
$r{-}s$ projections onto the direct
summands $A_l[l^k]$ having orders divisible by $l^M.$
\medskip

\noindent {\bf Step 3.} Pick $k$ big enough such that (3.4) holds
and such that there is a $\sigma$ as constructed in Step 2. Consider the diagram (3.5). We choose an element
$$\gamma \in G(
F_{l^{k+1}}(\frac{1}{l^{k}}{\hat\Lambda_{I}},\frac{1}{l^{k}}{\hat\Lambda_{J}})/F_{l^k})
 \subset G(F_{l^{k+1}}(\frac{1}{l^{k}}{\hat\Lambda_{I}},\frac{1}{l^{k}}{\hat\Lambda_{J}})/F),$$ 
in the following manner:
$\gamma |_{F_{l^k}(
\frac{1}{l^{k}}{\hat\Lambda_{I}},\frac{1}{l^{k}}{\hat\Lambda_{J}})}
= \sigma$ in the subgroup 
\newline
$G((F_{l^{k}}(\frac{1}{l^{k}}{\hat\Lambda_{I}},
\frac{1}{l^{k}}{\hat\Lambda_{J}})/
F_{l^{k}}(\frac{1}{l^{k}}{\hat\Lambda_{I}}))$ of 
$G( F_{l^{k}}(\frac{1}{l^{k}}{\hat\Lambda_{I}},
\frac{1}{l^{k}}{\hat\Lambda_{J}})/F_{l^k})$ and $\gamma
|_{F_{l^{k+1}}} = h$ is such that, the action of $h$ on the module
$T_l$ is given by a nontrivial homothety $1+l^ku_0,$ for some $u_0
\in \Z_{l}^{\times}.$ Such a homothety $h$ exists by the
assumption on the image  ${\rho}_{l}(G_F).$ By the Chebotarev
density theorem there are infinitely many prime ideals $v$ in
${\cal O}_F$ such that $\gamma$ is equal to the Frobenius element
for the prime $v$ in the extension
$F_{l^{k+1}}(\frac{1}{l^{k}}{\hat\Lambda_{I}},
\frac{1}{l^{k}}{\hat\Lambda_{J}})/F.$
In the remainder of the proof we work with prime ideals $v$ which
we have just selected. For each such $v$ we fix a prime $w$ in
$F_{l^{k+1}}(\frac{1}{l^{k}}{\hat\Lambda_{I}},
\frac{1}{l^{k}}{\hat\Lambda_{J}})$ above $v.$
\medskip

\noindent {\bf Step 4.} Let $i \in I$ and let $l^{c_i}$ be the
order of $r_v(P_i)$ in $B_v(\kappa_v)_l,$ for some $c_i \geq 0.$
Hence, $l^{k{+}c_i}\frac{1}{l^k}\, r_v( P_i) = 0.$ So $Q_i =
\frac{1}{l^k}P_i \in B(F_{l^{k+1}}(\frac{1}{l^k}\hat\Lambda_I))$
maps via the map $r_{w_1}$ (see the diagram (2.1)) to the point
$r_{w_1}(Q_i) \in B_v(\kappa_{w_1})_l$ of order $l^{k{+}c_i}.$
Here $w_1$ is a prime of $F_{l^{k+1}}(\frac{1}{l^k}\hat\Lambda_I)$
below $w.$  By the axioms ($A_2$), ($A_3$) and by the choice of
$v,$ the point $r_{w_1}(Q_i)$ comes from an element of
$B_v(\kappa_v)_l.$ Using the assumption that the right vertical
arrow in the diagram (2.1) is an isomorphism, by the choice of
$v,$ we see that the action of $h$ on $r_{w_1}(Q_i)$ is as
follows: $$h(r_{w_1}(Q_i)) = (1 + l^k u_0)r_{w_1}(Q_i).$$ But
$r_{w_1}(Q_i) \in B_v(\kappa_v)_l,$ so
$h(r_{w_1}(Q_i)){=}r_{w_1}(Q_i),$ again by the choice of $v.$
Hence, $l^k r_{w_1}(Q_i){=} 0.$ This can only happen, when $c_i =
0.$
\medskip

\noindent {\bf Step 5.} In this part of the proof we use the
argument similar to that in Proposition 2.19 of [BGK2]. Let $w_2$
denote the prime in
$F_{l^{k}}(\frac{1}{l^{k}}{\hat\Lambda_{I}},\frac{1}{l^{k}}{\hat\Lambda_{J}})$
below $w$ and let $u_2$ denote the prime in $F_{l^k}$ below $w_2.$
Consider the following commutative diagram: $$ \CD B(F)/l^kB(F)
@>>> B_v(\kappa_{v})/l^kB_v(\kappa_{v})\\ @VVV @VVV\\
B(F_{l^k})/l^kB(F_{l^k}) @>>>
B_v(\kappa_{u_2})/l^kB_v(\kappa_{u_2})\\ @VVV @VVV\\ 
Hom( (G_{F_{l^k}, S_{l}})^{ab}; A_l[l^k]) @>>> Hom (g_{u_2}; A_l[l^k]).
\endCD
\tag{3.8}$$
The bottom vertical maps are described in the natural
way using diagram 2.1 (similarly to the map $\Psi$ in the proof of Lemma 2.12). 
The bottom left vertical arrow is well defined since by axiom $A_1$
there is a natural isomorphism 
$$H^{1}_{f,S_{l}}(G_{L},T_{l}) \cong H^{1}(G_{L, S_{l}},T_{l}).$$
Every point
$P_j,$ where $j{\in} J,$ maps via the left vertical arrow in the diagram
(3.8) to an element ${{\phi}^{(k)}_{P_j}}$ defined by (2.3).
The homomorphism ${{\phi}^{(k)}_{P_j}}$ factors through the group $G(
F_{l^{k}}(\frac{1}{l^{k}}{\hat\Lambda_{I}},\frac{1}{l^{k}}{\hat\Lambda_{J}})/F_{l^k}).$ We
denote this factorization with the same symbol ${{\phi}^{(k)}_{P_j}}.$ By the
choice of the element $\gamma$ in Step 3 (see also the comments following the diagram (3.7))
we see that the element ${{\phi}^{(k)}_{P_j}}(\gamma
|_{F_{l^k}(\frac{1}{l^{k}}{\hat\Lambda_{I}},\frac{1}{l^{k}}{\hat\Lambda_{J}})})
\in A_l[l^k]$ has order divisible by $l^M.$ Hence,  the element
${{\phi}^{(k)}_{P_j}} \in Hom (H^{ab}_{l^k}; A_l[l^k])$ has order
divisible by $l^M,$ and by the choice of $v$ in Step 3, it maps via 
the bottom horizontal arrow in the diagram (3.8) to an element of order divisible by $l^M.$ So the point $P_j$ maps to an
element of order divisible by $l^M$ via the top horizontal arrow in the diagram (3.8). \qed
\enddemo

\medskip
\proclaim{Corollary 3.9} With the same assumptions as in the
Theorem 3.1 the Galois group $$
G(F_{l^{k}}(\frac{1}{l^{k}}\hat\Lambda_I) \cap
F_{l^{k}}(\frac{1}{l^{k}}\hat\Lambda_J)\,/\,F_{l^{k}})$$ has order
bounded independently of $k.$
\endproclaim

\demo{Proof}  In the same way as in Step 2 of the proof of Theorem 3.1
we show that the horizontal Kummer maps in the diagram (3.10) below
have cokernels bounded independently of $k$ (see diagram (3.7)).
$$ \diagram &
{\phantom{\Big|} G(F_{l^{k}}(\frac{1}{l^{k}}\hat\Lambda_I,
\frac{1}{l^{k}}\hat\Lambda_J)/F_{l^{k}}(\frac{1}{l^{k}}\hat\Lambda_J))\phantom{\Big|} }
 \ar@{^{(}->}[r]\ar@{^{(}->}[d]&{\phantom{\Big|} \bigoplus_{i=1}^sA_l[l^k]\phantom{\Big|} }\ar@{^{(}->}[d] &\\
 & {\phantom{\Big|} G(F_{l^{k}}(\frac{1}{l^{k}}\hat\Lambda_I,
\frac{1}{l^{k}}\hat\Lambda_J)/F_{l^{k}})\phantom{\Big|} }\ar@{^{(}->}[r] &
{\phantom{\Big|} \bigoplus_{i=1}^rA_l[l^k]\phantom{\Big|} }  &\\ &{\phantom{\Big|} G(F_{l^{k}}(\frac{1}{l^{k}}\hat\Lambda_I,
\frac{1}{l^{k}}\hat\Lambda_J)/F_{l^{k}}(\frac{1}{l^{k}}\hat\Lambda_I))\phantom{\Big|} }\ar@{^{(}->}[r]\ar@{^{(}->}[u]&
{\phantom{\Big|} \bigoplus_{i=1}^{r{-}s} A_l[l^k]\phantom{\Big|} } \ar@{^{(}->}[u]&\\
\enddiagram \tag{3.10}$$
The diagram (3.10) and the following diagram.
$$\diagram
&&&F_{l^{k}}(\frac{1}{l^{k}}{\hat\Lambda_{I}},\frac{1}{l^{k}}{\hat\Lambda_{J}})&
\\ &&F_{l^{k}}(\frac{1}{l^{k}}{\hat\Lambda_{I}})\urline&&
\ulline F_{l^{k}}(\frac{1}{l^{k}}{\hat\Lambda_{J}})
\\&&& F_{l^{k}}\ulline\urline&&\\
\enddiagram
$$ show that the subgroup $$
G(F_{l^{k}}(\frac{1}{l^{k}}\hat\Lambda_I,
\frac{1}{l^{k}}\hat\Lambda_J)/F_{l^{k}}(\frac{1}{l^{k}}\hat\Lambda_J))
\, G(F_{l^{k}}(\frac{1}{l^{k}}\hat\Lambda_I,
\frac{1}{l^{k}}\hat\Lambda_J)/F_{l^{k}}(\frac{1}{l^{k}}\hat\Lambda_I))$$
of $G(F_{l^{k}}(\frac{1}{l^{k}}\hat\Lambda_I,
\frac{1}{l^{k}}\hat\Lambda_J)/F_{l^{k}})$ has a finite index
bounded independently of $k.$ Hence, the Galois group
$$G(F_{l^{k}}(\frac{1}{l^{k}}\hat\Lambda_I) \cap
F_{l^{k}}(\frac{1}{l^{k}}\hat\Lambda_J)\,/\,F_{l^{k}})$$ has 
order bounded independently of $k.$ \qed
\enddemo

\proclaim{Lemma 3.11}
For any finite extension $L/F$ and  any prime
$w \notin S_l$ in ${\cal O}_{L},$ the reduction map
$$
\CD
r_{w}\,\,:\,\, B(L)_{tor}@>>> B_v(\kappa_w)$$
\endCD
$$
is an imbedding.
\endproclaim

\demo{Proof} Since ($A_1$) and ($A_2$) are fulfilled by
assumption, the lemma follows by the diagram (2.1) and  [BGK2],
Lemma 2.13. Note that Lemma 2.13 of [BGK2] holds also for
$l{=}2.\qed$
\enddemo

Under the assumptions of Section 2, we can derive the following stronger
version of Theorem 2.9 in the case when ${\cal O} = \Z.$

\proclaim{Theorem 3.12}

\noindent
Let $P$ and $P_1,\dots, P_r $ be nontorsion elements of $B(F)$ such
that $P_1,P_2,\dots, P_r$ are
linearly independent over ${\Z}.$ Denote by $\Lambda$ the submodule of 
$B(F)$ generated
by $P_1,P_2,\dots, P_r.$ The following two statements are equivalent:
\roster
\item[1]\,\,$P\in \Lambda$
\item[2]\,\,$r_{v}(P) \in r_{v}(\Lambda)$ for almost all primes $v$ of $F.$
\endroster
\endproclaim

\demo{Proof} We prove that $(2)$ implies $(1).$ The opposite
implication is obvious. By Theorem 2.9 there exists $a\in \Z$ such
that 
$$aP=\alpha_{1}P_{1}+\dots +\alpha_{r}P_{r},\tag{3.13}$$ 
for some $\alpha_i\in {\Z}.$ Let $l^M$ be the largest power of $l$
that divides $a.$ Fix $1 \leq i \leq r.$ By Theorem 3.1 there
are infinitely many primes $v$ such that $$r_v(P_1) = \dots =
r_v(P_{i-1}) = r_v(P_{i+1}) = \dots = r_v(P_r) = 0$$ in
$B_v(\kappa_v)_l$ and $r_v(P_i)$ has order at least $l^M$ in
$B(k_v)_l.$ Applying $r_v$ to the equality (3.13) we get: $$a
r_v(P) = \alpha_i r_v(P_i).\tag{3.14}$$ By the assumption of the
theorem and the choice of $v,$ 
$$r_v(P) = \beta_i r_v(P_i)\tag{3.15}$$ in $B_v(\kappa_v)_l,$ for some $\beta_i \in
\Z.$ Multipying both sides of (3.15) by $a$ and comparing with
(3.14) we get: $$(\alpha_i{-}a \beta_i) r_v(P_i) = 0$$ in
$B_v(\kappa_v)_l,$ which implies that $l^M$ divides $\alpha_i.$ So
the equality (3.13) gives: $${a \over l^M}P={\alpha_{1} \over
l^M}P_{1}+\dots +{\alpha_{r} \over l^M}P_{r} + R,\tag{3.16}$$
where $R \in B(F)[l^M].$ By Theorem 3.1 we pick infinitely
many primes $v$ in $F$ such that $r_v(P_1) {=} \dots {=} r_v(P_r)
{=} 0$ in $B_v(\kappa_v)_l,$ and apply the assumption that
$r_{v}(P) \in r_{v}(\Lambda),$ for almost all $v.$ Hence, by
applying $r_v$ to the equality (3.16), we get $r_v(R) {=} 0,$ for
infinitely many primes $v.$ This contradicts Lemma 3.11, unless
$R{=} 0.$ Hence, we obtain the equality: $${a \over
l^M}P={\alpha_{1} \over l^M}P_{1}+\dots+ {\alpha_{r} \over
l^M}P_{r}.$$ Repeating  the above argument for primes dividing ${a
\over l^M}$ we finish the proof by induction. \qed\enddemo
\bigskip

\subhead{4. Corollaries of Theorem 3.12}
\endsubhead
\bigskip

In the case of K-groups, we put 
$B(F){=}K_{2n+1}(F)/C_{F},$ for 
$n{\geq}0,$ where $C_F$ is the subgroup of 
$K_{2n+1}(F)$ generated by $l$-parts (for all primes $l$) 
of kernels of the Dwyer-Friedlander map. Observe that if the 
Quillen-Lichtenbaum conjecture holds true, then $B(F){=}K_{2n{+}1}(F)$ 
up to $2$-torsion. We obtain the following specialization of Theorem 3.12 
which strengthens the Theorem of [BGK1]. 
\bigskip

\proclaim{Theorem 4.1}

\noindent 
Let $P$ and $P_1,$ $P_2,$ $\dots,$ $P_r$
be nontorsion elements of $K_{2n+1}(F)/C_{F},$ such that 
$P_1,$ $P_2,$ $\dots,$ $P_r$ are
linearly independent over ${\Z}.$ Let $\Lambda\subset K_{2n+1}(F)/C_{F}$ 
denote the subgroup generated by $P_1,P_2,\dots, P_r.$ The following two 
statements are equivalent: 
\roster 
\item[1]\,\,$P\in \Lambda$ 
\item[2]\,\,$r_{v}(P) \in r_{v}(\Lambda)$ for almost all primes $v$ of $F.$
\endroster
\endproclaim

\demo{Proof}
Similarly to [BGK2], Example 3.4, let:
$$ \rho_l\,\,:\,\,G_F\rightarrow GL(T_l)\cong \Z_l^{\times}$$
be the one dimensional representation given by the $(n{+}1)$th
tensor power of the cyclotomic character. For every finite
extension $L/F,$ let $C_L$ be the subgroup of $K_{2n+1}(L)$ generated by
the $l$-parts (for all primes $l$) of the kernels of
the Dwyer-Friedlander maps cf. [DF]:
$$K_{2n+1}(L)\longrightarrow K_{2n+1}(L) \otimes \Bbb Z_l\longrightarrow
H^1(G_L;\, \Bbb Z_l(n+1)).$$ Let us put:
$$B(L)=K_{2n+1}(L)/C_L$$
and $c:=\#C_F.$ Because $T_l{=}\Z_l(n{+}1)$ is one{-}dimensional over $\Z_l$ and $B(F)$
is finitely generated, the axioms ($A_1){-}(A_3$) and
($B_1){-}(B_4$) are clearly satisfied (cf. [BGK2], Section 6).
In this case, we have $\cal O=\Z,$ ${\Cal G}^{alg}_{l}=\G_m$ and the map
${\cal G}_{l}^{alg}(\Z_{l})\rightarrow {\cal
G}_{l}^{alg}(\Z/l^{k})$ is the natural projection $\Z_{l}^{\times}\rightarrow ({\Z/l^{k}})^{\times}.$
The image of $\rho_l$
is open in ${\Z_{l}}^{\times},$ for each $l,$
which is clear from the definition of $\rho_l.
$\qed\enddemo

\proclaim{Theorem 4.2}

\noindent
Let $A$ be a principally polarized abelian variety of dimension
$g$ defined over the number field $F$ such that $End(A) = \Z$
and $dim(A)=g$ is either odd or $g=2$ or $6.$
Let $P$ and $P_1,\dots, P_r $ be nontorsion elements of $A(F)$ such that $P_1,P_2,\dots, P_r$ are
linearly independent over ${\Z}.$ Denote by $\Lambda$
the subgroup of $A(F)$ generated
by $P_1,P_2,\dots, P_r.$ Then the following two statements are equivalent:
\roster
\item[1]\,\,$P\in \Lambda$
\item[2]\,\,$r_{v}(P) \in r_{v}(\Lambda)$ for almost all primes $v$ of $F.$
\endroster
\endproclaim

\demo{Proof} Let $\rho_l:\,G_F\longrightarrow GL(T_l(A))$ be the
$l$-adic representation associated to $A.$ For a finite extension
$L/F$ let  $B(L){=}A(L).$ The axioms ($A_1){-}(A_3)$ are satisfied
by [BGK2], Examples 3.6, 3.7. The axioms ($B_1$) and ($B_2$) are
satisfied by the results of Faltings [F], Satz 4, and Zarhin [Z],
Corollary 5.4.5. The condition ($B_3$) holds due to the result of
Serre, [Se2], Corollary of Theorem 2, and the condition ($B_4$)
holds, because $B(F){=}A(F)$ is finitely generated by
theorem of Mordell and Weil. In this case,
we have ${\cal O}{=}\Z,$ ${\cal G}^{alg}_{l}{=}GSp_{2g}$ and the
map ${\cal G}_{l}^{alg}(\Z_{l})\rightarrow {\cal
G}_{l}^{alg}(\Z/l^{k})$ is the natural map $GSp_{2g}(\Z_{l})
\rightarrow GSp_{2g}({\Z/l^{k}}),$ which is surjective, for every
$k>0,$ by [A], Lemma 3.3.2 (1), p.135. The image of $\rho_l$ is
open in $GSp_{2g}(\Z_{l}),$ for each $l,$ by [Se1], Th\' eor\` eme
3, p. 97 and it contains an open subgroup of homotheties by the
theorem of Bogomolov, cf. [Bo], Corollary 1, p.702. \qed\enddemo
\medskip

\proclaim{Corollary 4.3} Let $E$ be an elliptic curve without complex multiplication, which is
defined over the number field $F.$ Let $P$ and $P_1,\dots, P_r $
be non-torsion elements of
the Mordell-Weil group $E(F)$ and such that $P_1,\dots, P_r $ are linearly independent over $\Z.$
Let $\Lambda $ be the subgroup generated by $P_1,\dots, P_r.$
Then the following two statements are equivalent:
\roster
\item[1]\,\,$P\in \Lambda$
\item[2]\,\,$r_{v}(P) \in r_{v}(\Lambda)$ for almost all primes $v$ of $F.$
\endroster
\endproclaim
\bigskip

\noindent
{\it Acknowledgements}:\quad
The authors thank the Isaac Newton Institute in Cambridge for the financial support during 
the visit in September 2002. 
The second author thanks the Max-Planck-Institut f{\" u}r 
Mathematik in Bonn for the hospitality and the financial support during the visit in 2002/2003. 
The research was partially
financed by a KBN grant.


\bigskip\bigskip

\Refs
\widestnumber\key{AAAA}

\ref\key A \by  A.N. Andrianov
\book Quadratic forms and Hecke operators
\publ Springer-Verlag
\vol Grundlehren der math. Wissenschaften 286
\yr 1987
\endref

\ref\key BGK1
\by G. Banaszak, W. Gajda, P. Kraso{\' n}
\paper A support problem for $K$-theory of number fields
\jour C. R. Acad. Sci. Paris S{\' e}r. 1 Math.
\vol 331 no. 3
\yr 2000
\pages 185-190
\endref

\ref\key BGK2 \by G. Banaszak, W. Gajda, P. Kraso{\' n} 
\paper Support problem for the intermediate Jacobians of $l$-adic representations. 
\jour Journal of Number Theory
\vol 100  no. 1,
\yr (2003) 
\pages 133--168.
\endref

\ref\key BK 
\by S. Bloch, K. Kato 
\paper L-functions and Tamagawa
numbers of motives,\, The Grothendieck Festschrift, 
\vol I 
\yr 1990 
\pages 333-400
\endref

\ref\key B
\by A. Borel
\paper Cohomologie de $SL_n$ et values de fonctiones zeta
\jour Ann. Acad. Scuola Normale Superiore
\yr 1974
\pages 613-636
\vol 7
\endref

\ref\key Bo
\by F.A. Bogomolov
\book Sur l'alg\' ebricit\' e des repr\' esentations $l$-adiques
\publ C.R.Acad.Sci. Paris S\' er. A-B
\vol 290
\pages A701-A703
\yr 1980
\endref

\ref\key C-RS \by C. Corralez-Rodrig{\'a}{\~ n}ez, R. Schoof
\paper Support problem and its elliptic analogue \jour Journal of
Number Theory \vol 64 \yr 1997 \pages 276-290
\endref

\ref\key DF \by W. Dwyer, E. Friedlander
\paper Algebraic and \'
etale K-theory
\jour Trans. Amer. Math. Soc.
\vol 292
\yr 1985
\pages  247-280
\endref

\ref\key Fa
\by G. Faltings
\paper Endlichkeitss\" atze f\" ur
abelsche Variet\" aten \" uber Zahlk\" orpern
\jour Inv. Math.
\vol 73
\yr 1983
\pages 349-366
\endref

\ref\key G-SOT \by I. Garcia-Selfa, M.A. Olalla, J.M. Tornero
\paper Computing the rational torsion of an elliptic curve using
Tate normal form
 \jour Journal of Number Theory
 \vol 96
 \yr 2002
 \pages 76-88
\endref


\ref\key J \by G. Janusz \book Algebraic number theory \publ
Academic Press, London and New York \yr 1973
\endref

\ref\key K \by Ch. Khare \paper Compatible systems of mod p Galois
representations \jour preprint
\endref


\ref\key KP \by Ch. Khare, D. Prasad \paper Reduction of abstract
homomorphisms of lattices mod $p$ \jour preprint \yr 2002 \pages
\endref

\ref\key L \by M. Larsen  \paper The support problem for abelian
varieties \jour preprint \yr 2002
\endref

\ref\key Q1 \by D. Quillen \paper Finite generation of the groups
$K_{i}$ of rings of algebraic integers. \jour Lecture Notes in
Mathematics \yr 1973 \vol 341 \pages 179-198
\endref

\ref
\key Q2
\by D. Quillen
\paper On the cohomology and K-theory of the general linear
group over finite field
\jour Ann. of Math.
\vol 96
\yr 1972
\pages 552-586
\endref

 \ref\key Ri \by K.A. Ribet \paper Kummer theory on extensions
of abelian varieties by tori \jour Duke Mathematical Journal \yr
1979 \vol 46, No. 4 \pages 745-761
\endref


\ref
\key RS \by K. Rubin, A. Silverberg
\paper Ranks of elliptic curves
\jour Bull. of the AMS
\yr  2002
\vol 39
\pages 455-474
\endref

\ref\key Se1
\by J.-P. Serre
\paper R\'esum\'es des cours au Coll\`ege de France
\jour Annuaire du Coll\`ege de France
\yr 1985-1986
\pages 95-99
\endref

\ref\key Se2
\by J.-P. Serre
\paper Sur les groupes de congruence des vari\'et\'es ab\'eliennes. II
\jour Izv. Akad. Nauk SSSR Ser. Mat.
\yr 1971
\vol 35
\pages 731-737
\endref


\ref\key Si \by J. Silverman \book The arithmetic of elliptic
curves \publ GTM 106, Springer-Verlag \yr 1986
\endref

\ref\key We \by T. Weston \paper Kummer theory of abelian
varieties and reductions of  Mordell-Weil groups \jour preprint
\yr 2002
\endref

\ref\key Za \by J.G. Zarhin \paper A finiteness theorem for
unpolarized Abelian varieties over number fields with prescribed
places of bad reduction \jour Invent. math. \yr 1985 \vol 79
\pages 309-321
\endref
\endRefs
\enddocument